\documentclass{amsart}
\usepackage{amssymb,amsmath,theorem,euscript}

\def\sm{\smallskip}

\newcounter{punct}

\def\punct{\refstepcounter{punct}{\arabic{punct}.  }}

\begin{document}

 \def\ov{\overline}
\def\wt{\widetilde}
 \newcommand{\rk}{\mathop {\mathrm {rk}}\nolimits}
\newcommand{\Aut}{\mathop {\mathrm {Aut}}\nolimits}
\newcommand{\Out}{\mathop {\mathrm {Out}}\nolimits}
\renewcommand{\Re}{\mathop {\mathrm {Re}}\nolimits}
\def\Br{\mathrm {Br}}

\def\SL{\mathrm {SL}}
\def\SU{\mathrm {SU}}
\def\GL{\mathrm {GL}}
\def\U{\mathrm U}
\def\OO{\mathrm O}
 \def\Sp{\mathrm {Sp}}
 \def\SO{\mathrm {SO}}
\def\SOS{\mathrm {SO}^*}
 \def\Diff{\mathrm{Diff}}
 \def\Vect{\mathfrak{Vect}}
\def\PGL{\mathrm {PGL}}
\def\PU{\mathrm {PU}}
\def\PSL{\mathrm {PSL}}
\def\Symp{\mathrm{Symp}}
\def\End{\mathrm{End}}
\def\Mor{\mathrm{Mor}}
\def\Aut{\mathrm{Aut}}
 \def\PB{\mathrm{PB}}
 \def\cA{\mathcal A}
\def\cB{\mathcal B}
\def\cC{\mathcal C}
\def\cD{\mathcal D}
\def\cE{\mathcal E}
\def\cF{\mathcal F}
\def\cG{\mathcal G}
\def\cH{\mathcal H}
\def\cJ{\mathcal J}
\def\cI{\mathcal I}
\def\cK{\mathcal K}
 \def\cL{\mathcal L}
\def\cM{\mathcal M}
\def\cN{\mathcal N}
 \def\cO{\mathcal O}
\def\cP{\mathcal P}
\def\cQ{\mathcal Q}
\def\cR{\mathcal R}
\def\cS{\mathcal S}
\def\cT{\mathcal T}
\def\cU{\mathcal U}
\def\cV{\mathcal V}
 \def\cW{\mathcal W}
\def\cX{\mathcal X}
 \def\cY{\mathcal Y}
 \def\cZ{\mathcal Z}
\def\0{{\ov 0}}
 \def\1{{\ov 1}}
 \def\frA{\mathfrak A}
 \def\frB{\mathfrak B}
\def\frC{\mathfrak C}
\def\frD{\mathfrak D}
\def\frE{\mathfrak E}
\def\frF{\mathfrak F}
\def\frG{\mathfrak G}
\def\frH{\mathfrak H}
\def\frI{\mathfrak I}
 \def\frJ{\mathfrak J}
 \def\frK{\mathfrak K}
 \def\frL{\mathfrak L}
\def\frM{\mathfrak M}
 \def\frN{\mathfrak N} \def\frO{\mathfrak O} \def\frP{\mathfrak P} \def\frQ{\mathfrak Q} \def\frR{\mathfrak R}
 \def\frS{\mathfrak S} \def\frT{\mathfrak T} \def\frU{\mathfrak U} \def\frV{\mathfrak V} \def\frW{\mathfrak W}
 \def\frX{\mathfrak X} \def\frY{\mathfrak Y} \def\frZ{\mathfrak Z} \def\fra{\mathfrak a} \def\frb{\mathfrak b}
 \def\frc{\mathfrak c} \def\frd{\mathfrak d} \def\fre{\mathfrak e} \def\frf{\mathfrak f} \def\frg{\mathfrak g}
 \def\frh{\mathfrak h} \def\fri{\mathfrak i} \def\frj{\mathfrak j} \def\frk{\mathfrak k} \def\frl{\mathfrak l}
 \def\frm{\mathfrak m} \def\frn{\mathfrak n} \def\fro{\mathfrak o} \def\frp{\mathfrak p} \def\frq{\mathfrak q}
 \def\frr{\mathfrak r} \def\frs{\mathfrak s} \def\frt{\mathfrak t} \def\fru{\mathfrak u} \def\frv{\mathfrak v}
 \def\frw{\mathfrak w} \def\frx{\mathfrak x} \def\fry{\mathfrak y} \def\frz{\mathfrak z} \def\frsp{\mathfrak{sp}}
 \def\bfa{\mathbf a} \def\bfb{\mathbf b} \def\bfc{\mathbf c} \def\bfd{\mathbf d} \def\bfe{\mathbf e} \def\bff{\mathbf f}
 \def\bfg{\mathbf g} \def\bfh{\mathbf h} \def\bfi{\mathbf i} \def\bfj{\mathbf j} \def\bfk{\mathbf k} \def\bfl{\mathbf l}
 \def\bfm{\mathbf m} \def\bfn{\mathbf n} \def\bfo{\mathbf o} \def\bfp{\mathbf p} \def\bfq{\mathbf q} \def\bfr{\mathbf r}
 \def\bfs{\mathbf s} \def\bft{\mathbf t} \def\bfu{\mathbf u} \def\bfv{\mathbf v} \def\bfw{\mathbf w} \def\bfx{\mathbf x}
 \def\bfy{\mathbf y} \def\bfz{\mathbf z} \def\bfA{\mathbf A} \def\bfB{\mathbf B} \def\bfC{\mathbf C} \def\bfD{\mathbf D}
 \def\bfE{\mathbf E} \def\bfF{\mathbf F} \def\bfG{\mathbf G} \def\bfH{\mathbf H} \def\bfI{\mathbf I} \def\bfJ{\mathbf J}
 \def\bfK{\mathbf K} \def\bfL{\mathbf L} \def\bfM{\mathbf M} \def\bfN{\mathbf N} \def\bfO{\mathbf O} \def\bfP{\mathbf P}
 \def\bfQ{\mathbf Q} \def\bfR{\mathbf R} \def\bfS{\mathbf S} \def\bfT{\mathbf T} \def\bfU{\mathbf U} \def\bfV{\mathbf V}
 \def\bfW{\mathbf W} \def\bfX{\mathbf X} \def\bfY{\mathbf Y} \def\bfZ{\mathbf Z} \def\bfw{\mathbf w}
 \def\R {{\mathbb R }} \def\C {{\mathbb C }} \def\Z{{\mathbb Z}} \def\H{{\mathbb H}} \def\K{{\mathbb K}}
 \def\N{{\mathbb N}} \def\Q{{\mathbb Q}} \def\A{{\mathbb A}} \def\T{\mathbb T} \def\P{\mathbb P} \def\G{\mathbb G}
 \def\bbA{\mathbb A} \def\bbB{\mathbb B} \def\bbD{\mathbb D} \def\bbE{\mathbb E} \def\bbF{\mathbb F} \def\bbG{\mathbb G}
 \def\bbI{\mathbb I} \def\bbJ{\mathbb J} \def\bbL{\mathbb L} \def\bbM{\mathbb M} \def\bbN{\mathbb N} \def\bbO{\mathbb O}
 \def\bbP{\mathbb P} \def\bbQ{\mathbb Q} \def\bbS{\mathbb S} \def\bbT{\mathbb T} \def\bbU{\mathbb U} \def\bbV{\mathbb V}
 \def\bbW{\mathbb W} \def\bbX{\mathbb X} \def\bbY{\mathbb Y} \def\kappa{\varkappa} \def\epsilon{\varepsilon}
 \def\phi{\varphi} \def\le{\leqslant} \def\ge{\geqslant}

\def\UU{\bbU}
\def\Mat{\mathrm{Mat}}
\def\tto{\rightrightarrows}

\def\Gr{\mathrm{Gr}}
\def\Isom{\mathrm{Isom}}

\def\graph{\mathrm{graph}}

\def\O{\mathrm{O}}

\def\la{\langle}
\def\ra{\rangle}

\def\B{\mathrm B}
\def\Int{\mathrm {Int}}

\begin{center}
\Large\bf
A remark on representations of infinite symmetric groups
\bigskip

\large\sc
Yu.A.Neretin\footnote{
\newline
Supported by the grant FWF, P22122.}
\end{center}

{\small We simplify  construction of Thoma representations of an infinite symmetric group.}

\bigskip

{\bf\punct Spherical representations.} Let $G$ be a group, $K$ a subgroup.
The pair $(G,K)$ is called  {\it spherical} if for any irreducible unitary representation
$\rho$
of $G$   the dimension of the subspace of $K$-fixed vectors is $\le 1$. An unit vector
of this subspace is called a {\it spherical vector.}  Recall that the
{\it spherical function} of an irreducible spherical representation $\rho$ of $G$ is given by
$$
\Phi(g):=\la \rho(g) v,v\ra
,$$
 where $v$ is the unit $K$-fixed vector.

\sm

{\bf\punct The double of symmetric group.} Let $\Omega$ be a countable set.  A permutation of
$\Omega$
is called {\it finite} if it fixes all but finite number elements of $\Omega$
 Denote by $S(\Omega)$ the group of finite permutations of a countable set $\Omega$.
 We also use an alternative notation $S_\infty$ for such groups.
 
Now let $G=S_\infty\times S_\infty$ be the product of two copies of $S_\infty$,  
let $K\simeq S_\infty$ be the diagonal subgroup.
By \cite{Olsh}, the pair $(G,K)$ is spherical.

\sm

{\bf\punct Thoma formula.} Let $G$, $K$ be the same as in 
the previous subsection.
  Representations of $G$ spherical with respect to $K$
 are parametrized by collection of positive numbers
$$
\alpha_1\ge\alpha_2\ge\dots,\qquad \beta_1\ge \beta_2\ge \dots, \qquad \sum \alpha_i+\sum\beta_j\le1
,$$
finite and empty collections are admissible.
Spherical functions  are given by the formula
\begin{equation}
\Phi_{\alpha,\beta}(\sigma,\tau)=
\prod_{k=2} \left(\sum_i \alpha_i^k+(-1)^{k-1}\sum_j
 \beta_j^k\right)^{\text{number of cycles of $\sigma\tau^{-1}$ of length $k$}}
 \label{eq:thoma}
\end{equation}
(the product is finite),
 Thoma (1964) formulated this statement on another  language (see \cite{Tho}, \cite{Olsh}).

\sm

Explicit constructions of Thoma representations 
 were obtained by Vershik and Kerov (see \cite{VK}, another version was done by Vershik \cite{Ver}).
 Olshanski \cite{Olsh} proposed a transparent construction (which also is a modification 
 of \cite{VK}), but it does not cover all possible values of parameters, we 
 add a missing element to his construction.

{\bf\punct Olshanski construction.}
Let
\begin{equation}
\sum \alpha_i+\sum\beta_j=1
\label{eq:=1}
\end{equation}
Consider a Hilbert space $H=H_{\ov 0}\oplus H_{\ov 1}$. Fix an orthonormal basis
$e_i$ in $H_{\ov 0}$ and $f_j\in H_{\ov 1}$. Consider a unit vector 
$\xi\in H\otimes H$ given by
$$
\xi:=
\sum_i \alpha_i^{1/2} e_i\otimes e_i+ \sum_{j} \beta_j^{1/2} f_j\otimes f_j.
$$
We consider the infinite super-tensor product%
\footnote{Consider a $\Z_2$-graded space
 $V=V_{\ov 0}\oplus V_{\ov 1}$. Then $V\otimes V$ is the usual tensor product, but the operator
 of transposition of factors is another: for homogeneous elements $v$, $w$ we have $v\otimes w\to (-1)^{\epsilon} w\otimes v$, where
 $\epsilon=1$ if $v$, $w\in V_{\ov 1}$ and $\epsilon=0$ if at least one of the vectors $v$, $w$
 is contained in $V_{\ov 0}$. On $n$-factor product $V\otimes\dots\otimes V$ we have an action 
 of the symmetric group $S_n$, any permutation is a product of transpositions, and action of
 a transposition was described above.}
\begin{equation}
\cH:=(H\otimes H,\xi)\otimes (H\otimes H,\xi)\otimes (H\otimes H,\xi)\otimes\dots
\label{eq:cH}
\end{equation}
(recall that the definition of an infinite tensor product
 requires distinguished unit vectors). The subgroup 
$S_\infty\times e\subset G$ acts in $\cH$ by permutations of first factors in brackets 
$(H\otimes H,\xi)$, the group $e\times S_\infty$ acts by permutations of the second factors.
The diagonal subgroup $K$ acts by permutations of the whole brackets $(H\otimes H,\xi)$
and $\xi^{\otimes \infty}$ is the $K$-spherical vector.

However, the condition (\ref{eq:=1}) is essential%
\footnote{Otherwise the length of $\xi$ is not 1 and (\ref{eq:cH})
is not well-defined.}, the same dif\-fi\-cul\-ty arises for other 
spherical pairs considered by Olshanski \cite{Olsh} and in a more general 
setting discussed in \cite{Ner1}, \cite{Ner2}.

\sm

{\bf\punct  Products of spherical functions.}
 Consider two irreducible spherical representations of $G$ corresponding to
parameters $\alpha$, $\beta$ and $\alpha'$, $\beta'$. It can be easily shown that their tensor product
contains only one $K$-fixed vector%
\footnote{but this representation can be reducible}, and the spherical function is the product of spherical functions,
it has the  form (\ref{eq:thoma}) with parameters
$$
\wt\alpha=\{\alpha_i \alpha'_k\},\,\{\beta_j\beta'_l\},
\qquad
\wt\beta=\{\alpha_i\beta_l'\},\,\{\beta_j\alpha_k'\}
.$$
Therefore {\it to construct all spherical representations it suffices to construct a representation%
\footnote{Another construction of this representation can be found in \cite{DN}.}
corresponding to a single-element collection $\alpha$ and empty collection $\beta$},
the spherical function of this representation is
\begin{equation}
\Psi(\sigma,\tau)=\alpha^{\text{number of $i$ such that $\sigma i\ne\tau i$ }}
.
\label{eq:Psi}
\end{equation}

\sm

{\bf\punct Group of affine isometries and Araki scheme.} For details, see, e.g., \cite{Ner-book}, V.1.6-7, X.1.
 Denote by $F_n$ the Hilbert space of holomorphic
functions on $\C^n$ with the inner product 
$$
\la f, g\ra=\frac 1{\pi^n}\int_{\C^n} f(z)\ov{g(z)}\,e^{-\la z,z\ra}dz\,d\ov z
.
$$
Consider the natural embeddings $J_n:F_n\to F_{n+1}$ given by
$$
J_n f(z_1,\dots,z_n, z_{n+1})=f(z_1,\dots,z_n)
.$$
Evidently, $J_n$ is an isometric embedding. We consider the union of the chain
$$...\longrightarrow F_n \longrightarrow F_{n+1}\longrightarrow\dots$$
and its completion $\bfF_\infty$ (the {\it boson Fock space}), 
see, e.g., \cite{Ner-book}, VI.1.

Consider a real Hilbert space $H$. Denote by $\O(H)$ the group of orthogonal
(i.e., real unitary) operators.  Consider the group 
$\Isom(H)$
generated by $\O(H)$ and translations, we get a semi-direct product
$\Isom(H)=\O(H)\ltimes H$, this group acts in $H$ by affine transformations
$$
h\mapsto Ah+v, \qquad \text{where $A\in\O(H)$, $v\in H$}
.
$$

Next, we construct the linear representation $Exp\,(\cdot)$ of $\Isom(\ell_2)$ in the Fock space $\bfF_\infty$.
Orthogonal transformations act in the Fock space by
$$
Exp\,(A)f(z)= f(zA), \qquad A\in\O(\ell_2)
,$$
translations by
$$
Exp\,(v)f(z)=f(z+v) e^{-\la z,v\ra-\frac 12\la v,v\ra},\qquad v\in\ell_2
.
$$
Thus we get an irreducible  unitary representation of $\Isom(\ell_2)$. The function $f(z)=1$ is  $\O(\ell_2)$-invariant,
the spherical function is $e^{-\frac 12\la h,h\ra}$.

One of the most common ways%
\footnote{See numerous examples in \cite{Ner-book}, Sections VIII.6.8, IX.1.6, IX.2.5, IX.5.4, IX.1.4, IX.1.5, IX.3.12, IX.4.6, F.4.}
 ({\it Araki scheme}) 
to construct representations of infinite-dimensional groups $G$
is embeddings of $G$ to the group of isometries of Hilbert  space and restrictions of
the representation 
$Exp(\cdot)$ to $G$. 

Let we have a unitary representation $U$ of a group $G$ in
a Hilbert space $H$. Let $\Xi:G\to H$ be a function  satisfying
\begin{equation}
\Xi(g_1g_2)=T(g_1)\Xi(g_2)+\Xi(g_1)
\label{eq:cocycle}
.\end{equation}
 Then affine isometric transformations
 \begin{equation}
 \wt U (g)h=U(g) h+\Xi(g)
 \label{eq:wtU}
 \end{equation}
satisfy 
$$\wt U(g_1)\wt U(g_2)=\wt U(g_1g_2)$$
 and we get embedding 
$G\to\Isom(H)$ (this is straightforward). For a fixed vector $\eta\in H$
the function
\begin{equation}
\Xi(g):=U(g)\eta-\eta
\label{eq:shift}
\end{equation}
satisfies the equation (\ref{eq:cocycle}). This solution of (\ref{eq:cocycle})
  is not interesting because such correction
$\Xi(\cdot)$ is equivalent to a change of origin of coordinates.

Now let $G$ acts by linear transformations of a larger linear space $\widehat H\supset H$, 
$\eta\in \widehat H\setminus H$. If $U(g)\eta-\eta\in H$ for all $g$, 
then we get a nontrivial affine isometric action of $G$.

\sm

{\bf\punct The construction for  the double of the symmetric  group.} 
Consider the space $\ell_2$ with basis $e_j$. The   group $G=S_\infty\times S_\infty$
 acts in $\ell_2\otimes\ell_2$ by linear transformations
$$
U(\sigma,\tau) e_i\otimes e_j=e_{\sigma i}\otimes e_{\tau j}
$$
We fix $s>0$, set
$$
\eta:=s\sum_{j=1}^\infty e_j\otimes e_j
$$
and define $\Xi(\sigma,\tau)$ by (\ref{eq:shift}).
Note that $\eta\notin \ell_2\otimes\ell_2$ but $\Xi(\sigma,\tau)\in\ell_2\otimes\ell_2$.
Also, $\Xi(\sigma,\sigma)=0$.
We define affine isometric action of $G= S_\infty\times S_\infty$ by 
(\ref{eq:wtU})
and restrict the representation of $\Isom(\ell_2\otimes\ell_2)$ to $G$.
Then the function $f=1$ is a unique $K$-fixed vector, the spherical function
is  (\ref{eq:Psi}) with $\alpha=e^{-s^2}$.

\sm

{\bf\punct Symmetric group and hyper-octahedral subgroup.}
Consider two copies of $\N$, say $\N_+$ and $\N_-$. Denote their points by
$1_+$, $2_+$, \dots, $1_-$, $2_-$, \dots. Let $G_1=S(\N_+\sqcup \N_-)=:S_{2\infty}$.
The {\it hyperoctohedral group} $K_1$ is the subgroup in $G_1$ consisting of permutations
such that for any $j\in \N$ the pair $(\sigma j_+, \sigma j_-)$ has the form
$(m_+, m_-)$ or $(m_-, m_+)$. Evidently, $K_1$ is a semidirect product,
$K_1=S_\infty\ltimes \Z_2^\infty$.

The pair $(G_1,K_1)$ is spherical, see \cite{Olsh}.

Now we consider the  Hilbert spaces $\ell_2(\N_\pm)$ with bases $e_j^\pm$. Consider the
Hilbert space
$$ 
H=\bigl(\ell_2(\N_+)\oplus \ell_2(\N_-)\bigr)\otimes \bigl(\ell_2(\N_+)\oplus \ell_2(\N_-)\bigr)
$$
The group $G_1=S(\N_+\sqcup \N_-)$ acts in this space in a natural way  (on each tensor factor). Next,
we define the vector 
$$
\eta:=  s\cdot\sum_j (e_j^+ \oplus e_j^- + e_j^- \oplus e_j^+)
$$
and construct an embedding of $G_1=S(\N_+\sqcup \N_-)$ to $\Isom(H)$ as above.

This gives $K_1$-spherical representations of $G_1$, which were not covered by explicit construction 
of \cite{Olsh}.

\sm

{\bf\punct One more example.} Consider the pair $G_2\supset K_2$, where $G_2=G_1$ and $K_2\subset K_1$
is the group of permutations $\sigma$ sending any ordered pair $(j_+,j_-)$ to a pair
$(m_+,m_-)$. In fact,
$$
K_2=S_\infty\subset S_\infty\ltimes \Z_2^\infty = K_1
$$
This pair is spherical (see \cite{Ner2}, note that this fact has not counterpart for finite symmetric groups). Now we fix real parameters $s$, $t$, set
$$
\eta=s\sum e_j^+\otimes e_j^- + t\sum e_j^-\otimes e_j^+
$$
and repeat the same arguments. 

\sm

{\bf\punct Triple products.} Now let $G_3=S_\infty\times S_\infty\times S_\infty$,
$K_3\simeq S_\infty$ be the diagonal. This pair is spherical, see \cite{Ner1}, \cite{Ner2}. We set
$$
H:=\ell_2\otimes\ell_2\otimes \ell_2
$$
and
$$
\eta= s\sum_j e_j\otimes e_j\otimes e_j.
$$

{\tt  Math.Dept., University of Vienna,

Institute for Theoretical and Experimental Physics, Moscow

Mech.Math. Dept., Moscow State University,

e-mail: neretin(at) mccme.ru

URL:www.mat.univie.ac.at/$\sim$neretin}

\end{document}